\documentclass[12pt]{amsart}

\usepackage[russian]{babel}
\usepackage{amsmath,amsfonts,amscd,amssymb,latexsym}
\font\teneufm=eufm10 \font\seveneufm=eufm7 \font\fiveeufm=eufm5
\newfam\frakturfam
\textfont\frakturfam=\teneufm \scriptfont\frakturfam=\seveneufm
\scriptscriptfont\frakturfam=\fiveeufm

\newtheorem{lemma}{╦хььр}
\newtheorem{theorem}{╥хюЁхьр}
\newtheorem{corol}{╤ыхфёЄтшх}

\def\bee{\begin{eqnarray}}
\def\bes{\begin{eqnarray*}}
\def\eee{\end{eqnarray}}
\def\ees{\end{eqnarray*}}

\def\a{\alpha}

\def\Proof{{\sl ─юърчрЄхы№ёЄтю.}\ }

\begin{document}

\pagestyle{plain}
\title{╨єўэ√х ш фшъшх ртЄюьюЁЇшчь√ рыухсЁ√ фшЇЇхЁхэЎшры№э√ї ьэюуюўыхэют Ёрэур 2}

\date{}

\maketitle

\begin{center}
{\bf ┴.└. ─єщёхэурышхтр$^{1}$, └.╤. ═рєЁрчсхъютр$^{1}$, ╙.╙. ╙ьшЁсрхт$^{2}$ }\\
{\normalsize{\sl $^{1}$ ┼тЁрчшщёъшщ эрЎшюэры№э√щ єэштхЁёшЄхЄ шь. ╦.═. ├єьшыхтр\\
 └ёЄрэр, 010008, ╩рчрїёЄрэ\\
e-mail: {\em bibinur.88@mail.ru}, {\em altyngul.82@mail.ru} }}

{\normalsize{\sl $^{2}$ Wayne State University\\
Detroit, MI 48202, USA\\
e-mail: {\em umirbaev@math.wayne.edu}}}
\end{center}

\begin{abstract}
─юърчрэю, ўЄю уЁєяяр Ёєўэ√ї ртЄюьюЁЇшчьют рыухсЁ√ фшЇЇхЁхэЎшры№э√ї ьэюуюўыхэют $k\{x,y\}$ эрф яюыхь $k$ їрЁръЄхЁшёЄшъш $0$ юЄ фтєї яхЁхьхээ√ї $x, y$ ё $m$ ъюььєЄшЁє■∙шьш фшЇЇхЁхэЎшЁютрэш ьш $\delta_1, \ldots, \delta_m$  ты хЄё  ётюсюфэ√ь яЁюшчтхфхэшхь ё юс·хфшэхэшхь. ╧юёЄЁюхэ яЁшьхЁ фшъюую ртЄюьюЁЇшчьр рыухсЁ√ $k\{x,y\}$ т ёыєўрх $m\geq 2$.
\end{abstract}

\noindent

{\bf ╩ы■ўхт√х ёыютр:} рыухсЁр фшЇЇхЁхэЎшры№э√ї ьэюуюўыхэют, Ёєўэ√х ш фшъшх ртЄюьюЁЇшчь√, ётюсюфэюх яЁюшчтхфхэшх.

\section{┬тхфхэшх}
\hspace*{\parindent}

╒юЁю°ю шчтхёЄэю \cite{Jung, Kulk, Czer, M-L1}, ўЄю ртЄюьюЁЇшчь√ рыухсЁ√ ьэюуюўыхэют $k[x,y]$ ш ётюсюфэющ рёёюЎшрЄштэющ рыухсЁ√ $k\left\langle x,y\right\rangle$ юЄ фтєї яхЁхьхээ√ї эрф яЁюшчтюы№э√ь яюыхь $k$  ты ■Єё  Ёєўэ√ьш.
┴юыхх Єюую \cite{Czer, M-L1}, уЁєяя√ ртЄюьюЁЇшчьют рыухсЁ  $k[x,y]$ ш $k\left\langle x,y\right\rangle$ шчюьюЁЇэ√, Є.х.
$$Aut_k k[x,y]\cong Aut_k k\left\langle x,y\right\rangle.$$

╚чтхёЄэю Єръцх, ўЄю ртЄюьюЁЇшчь√ фтєяюЁюцфхээ√ї ётюсюфэ√ї рыухсЁ ╧єрёёюэр эрф яюы ьш эєыхтющ їрЁръЄхЁшёЄшъш \cite{M-L2} ш ртЄюьюЁЇшчь√ фтєяюЁюцфхээ√ї ётюсюфэ√ї яЁртюёшььхЄЁшўэ√ї рыухсЁ эрф яЁюшчтюы№э√ьш яюы ьш \cite{Koz}  ты ■Єё  Ёєўэ√ьш. ╧. ╩юэ \cite{Cohn64} фюърчры, ўЄю ртЄюьюЁЇшчь√ ётюсюфэ√ї рыухсЁ ╦ш ъюэхўэюую Ёрэур  ты ■Єё  Ёєўэ√ьш. └эрыюу ¤Єющ ЄхюЁхь√ тхЁхэ фы  ётюсюфэ√ї рыухсЁ ы■сюую юфэюЁюфэюую °ЁрщхЁютюую ьэюуююсЁрчш  рыухсЁ \cite{Lewin}.
═ряюьэшь, ўЄю °ЁрщхЁют√ьш  ты ■Єё  ьэюуююсЁрчш  тёхї эхрёёюЎшрЄштэ√ї рыухсЁ \cite{14}, ъюььєЄрЄштэ√ї ш рэЄшъюььєЄрЄштэ√ї рыухсЁ \cite{15}, рыухсЁ ╦ш \cite{16, 17} ш ёєяхЁрыухсЁ ╦ш \cite{18,19}.

├Ёєяя√ ртЄюьюЁЇшчьют рыухсЁ ьэюуюўыхэют \cite{Shest1, Shest2, Shest3} ш ётюсюфэ√ї рёёюЎшрЄштэ√ї рыухсЁ \cite{Umir1, Umir2} юЄ ЄЁхї яхЁхьхээ√ї эрф яюыхь эєыхтющ їрЁръЄхЁшёЄшъш эх ьюуєЄ с√Є№ яюЁюцфхэ√ тёхьш ¤ыхьхэЄрЁэ√ьш ртЄюьюЁЇшчьрьш, Є.х. ёє∙хёЄтє■Є фшъшх ртЄюьюЁЇшчь√.  ╙.╙. ╙ьшЁсрхт√ь с√ыю фюърчрэю \cite{Umir1, Umir2}, ўЄю ртЄюьюЁЇшчь └эшър
$$\delta=(x+z(xz-zy), y+(xz-zy)z, z)$$
ётюсюфэющ рёёюЎшрЄштэющ рыухсЁ√ $k\left\langle x,y,z\right\rangle$ эрф яюыхь їрЁръЄхЁшёЄшъш $0$  ты хЄё  фшъшь.

╬ёэютэ√х яюэ Єш  фшЇЇхЁхэЎшры№э√ї рыухсЁ ьюцэю эрщЄш т ЁрсюЄрї \cite{Kol, KLMP, Ritt}.
╠√ сєфхь ЁрёёьрЄЁштрЄ№ фшЇЇхЁхэЎшры№э√х рыухсЁ√ ё ьэюцхёЄтюь ъюььєЄшЁє■∙шї фшЇЇхЁхэЎшЁютрэшщ $\Delta=\{\delta_1,\delta_2,\ldots,\delta_m\}$.
╧єёЄ№ $k$ --  фшЇЇхЁхэЎшры№эюх яюых їрЁръЄхЁшёЄшъш $0$ ш $k\{x,y\}$ -- рыухсЁр фшЇЇхЁхэЎшры№э√ї ьэюуюўыхэют эрф $k$ юЄ фтєї яхЁхьхээ√ї $x,y$.  ┼ёыш $|\Delta|=0$, Єю $k\{x,y\}$ ёЄрэютшЄё  юс√ўэющ рыухсЁющ ьэюуюўыхэют $k[x,y]$ эрф яюыхь $k$.
┬ ЁрсюЄрї ┬. трэ фхЁ ╩рыър \cite{Kulk} ш ╠. ═рурЄ√ \cite{Nag} фюърчрэю, ўЄю уЁєяяр $Aut(k[x,y])$  яЁхфёЄрты хЄё  т тшфх рьры№урьшЁютрээюую ётюсюфэюую яЁюшчтхфхэш 
$$Aut(k[x,y])=A \ast_{C} B,$$
уфх $A$ -- яюфуЁєяяр рЇЇшээ√ї ртЄюьюЁЇшчьют, $B$ -- яюфуЁєяяр ЄЁхєуюы№э√ї ртЄюьюЁЇшчьют ш $C=A\cap B$.

┬ эрёЄю ∙хщ ЁрсюЄх фюърч√трхЄё , ўЄю уЁєяяр Ёєўэ√ї ртЄюьюЁЇшчьют рыухсЁ√ $k\{x,y\}$ фюяєёърхЄ рэрыюушўэє■ ёЄЁєъЄєЁє рьры№урьшЁютрээюую ётюсюфэюую яЁюшчтхфхэш  фы  ы■сюую ьэюцхёЄтр фшЇЇхЁхэЎшЁютрэшщ $\Delta$.  ╩Ёюьх Єюую, шёяюы№чє  ¤Єє ёЄЁєъЄєЁє, фрхЄё  яЁшьхЁ фшъюую ртЄюьюЁЇшчьр рыухсЁ√ $k\{x,y\}$ яЁш $|\Delta|\geq 2$. ▌ЄюЄ яЁшьхЁ  ты хЄё  рэрыюуюь шчтхёЄэюую ртЄюьюЁЇшчьр └эшър \cite[ёЄЁ. 398]{Cohn}.

╥ръшь юсЁрчюь, ртЄюьюЁЇшчь√ рыухсЁ√ $k\{x,y\}$  ты ■Єё  Ёєўэ√ьш яЁш $|\Delta|=0$ ш $k\{x,y\}$ шьххЄ фшъшх ртЄюьюЁЇшчь√ яЁш $|\Delta|\geq 2$. ┬юяЁюё ю Ёєўэ√ї ш фшъшї ртЄюьюЁЇшчьрї рыухсЁ√ $k\{x,y\}$ юёЄрхЄё  юЄъЁ√Є√ь яЁш $|\Delta|=1$.

╤ЄрЄ№  юЁурэшчютрэр ёыхфє■∙шь юсЁрчюь. ┬ Ёрчфхых 2 яЁштхфхэ√ эхюсїюфшь√х юяЁхфхыхэш  ш ёЇюЁьєышЁютрэ√  эхъюЄюЁ√х шчтхёЄэ√х єЄтхЁцфхэш . ╨рчфхы 3 яюёт ∙хэ яЁхфёЄртыхэш■ уЁєяя√ Ёєўэ√ї ртЄюьюЁЇшчьют рыухсЁ√ $k\{x,y\}$ т тшфх рьры№урьшЁютрээюую ётюсюфэюую яЁюшчтхфхэш . ┬ Ёрчфхых 4 фюърч√трхЄё  ёюъЁрЄшьюёЄ№ ы■сюую эх рЇЇшээюую Ёєўэюую ртЄюьюЁЇшчьр рыухсЁ√ $k\{x,y\}$. ┬ Ёрчфхых 5 фрхЄё  яЁшьхЁ фшъюую ртЄюьюЁЇшчьр.

\section{╬яЁхфхыхэш  ш яЁхфтрЁшЄхы№э√х ЇръЄ√}
\hspace*{\parindent}

╧єёЄ№ $R$ -- яЁюшчтюы№эюх ъюььєЄрЄштэюх ъюы№Ўю ё хфшэшЎхщ. ╬ЄюсЁрцхэшх $d: R\rightarrow R$ эрч√трхЄё  \textit{фшЇЇхЁхэЎшЁютрэшхь}, хёыш фы  тёхї $s,t\in R$ т√яюыэ ■Єё  єёыютш 
$$d(s+t)=d(s)+d(t),$$
$$d(st)=d(s)t+sd(t).$$

╧єёЄ№ $\Delta=\{ \delta_1, \ldots, \delta_m \}$ -- юёэютэюх ьэюцхёЄтю фшЇЇхЁхэЎшры№э√ї юяхЁрЄюЁют.

╩юы№Ўю $R$ эрч√трхЄё  \textit{фшЇЇхЁхэЎшры№э√ь ъюы№Ўюь} шыш \textit{$\Delta$-ъюы№Ўюь}, хёыш $\delta_1, \ldots, \delta_m$   ты ■Єё  ъюььєЄшЁє■∙шьш фшЇЇхЁхэЎшЁютрэш ьш ъюы№Ўр $R$, Є.х.  $\delta_i: R\rightarrow R$ -- фшЇЇхЁхэЎшЁютрэш  ш $\delta_i \delta_j=\delta_j \delta_i$ фы  тёхї $i,j$.

╧єёЄ№ $\Theta$ -- ётюсюфэ√щ ъюььєЄрЄштэ√щ ьюэюшф эр ьэюцхёЄтх фшЇЇхЁхэЎшры№э√ї юяхЁрЄюЁют $\Delta=\{ \delta_1, \ldots, \delta_m \}$. ▌ыхьхэЄ√
$$\theta=\delta_1^{i_1}\ldots \delta_m^{i_m}$$
ьюэюшфр $\Theta$ эрч√тр■Єё  \textit{яЁюшчтюфэ√ьш юяхЁрЄюЁрьш}.  {\em ╧юЁ фъюь} $\theta$ эрч√трхЄё  ўшёыю  $|\theta|=i_1+\ldots+i_m$. ╧юыюцшь Єръцх $\gamma(\theta)=(i_1,\ldots,i_m)\in \mathbb{Z}_+^m$, уфх $\mathbb{Z}_+$ -- ьэюцхёЄтю тёхї эхюЄЁшЎрЄхы№э√ї Ўхы√ї ўшёхы.

╧єёЄ№ $R$ -- яЁюшчтюы№эюх фшЇЇхЁхэЎшры№эюх ъюы№Ўю ш яєёЄ№ $X=\{ x_1, \ldots, x_n\}$ --  ьэюцхёЄтю ёшьтюыют. ╨рёёьюЄЁшь ьэюцхёЄтю ёшьтюыют $X^\Theta=\{ x_i^\theta | 1\leq i \leq n, \theta\in \Theta\}$ ш рыухсЁє ьэюуюўыхэют $R[X^\Theta]$ эр ьэюцхёЄтх ёшьтюыют $X^\Theta$. ╧юырур  $$\delta_i(x_j^\theta)=x_j^{\theta \delta_i}$$
фы  тёхї $1\leq i \leq m, 1\leq j \leq n, \theta \in \Theta$, яЁхтЁрЄшь рыухсЁє $R[X^\Theta]$ т фшЇЇхЁхэЎшры№эє■ рыухсЁє. ─шЇЇхЁхэЎшры№эр  рыухсЁр $R[X^\Theta]$ юсючэрўрхЄё  ўхЁхч $R\{X\}$ ш эрч√трхЄё  \textit{рыухсЁющ фшЇЇхЁхэЎшры№э√ї  ьэюуюўыхэют} эрф $R$ юЄ ьэюцхёЄтр яхЁхьхээ√ї $X$ \cite{Kol}.

╧єёЄ№ $M$ -- ётюсюфэ√щ ъюььєЄрЄштэ√щ ьюэюшф юЄ ьэюцхёЄтр яхЁхьхээ√ї $x_i^{\theta}$, уфх $1\leq i\leq n$ ш $\theta\in \Theta$. ▌ыхьхэЄ√ $M$ эрч√тр■Єё   {\em ьюэюьрьш} рыухсЁ√  $R\{x_1,x_2,\ldots,x_n\}$. ╦■сющ ¤ыхьхэЄ
$a\in R\{x_1,x_2,\ldots,x_n\}$ юфэючэрўэю чряшё√трхЄё  т тшфх
\bes
a=\sum_{m\in M} r_m m
\ees
ё ъюэхўэ√ь ўшёыюь эхэєыхт√ї $r_m\in R$.

─ы  ы■сюую $x_i^{\theta}\in X^{\Theta}$ яюыюцшь  $\alpha(x_i^{\theta})=(\varepsilon_i,\gamma(\theta))\in \mathbb{Z}_+^{n+m}$, уфх $\varepsilon_1,\ldots,\varepsilon_n$ -- ёЄрэфрЁЄэ√щ срчшё $\mathbb{Z}_+^n$. ┼ёыш $m=a_1\ldots a_s\in M$, уфх $a_1,\ldots,a_s\in X^{\Theta}$, Єю яюыюцшь $\alpha(m)=\alpha(a_1)+\ldots +\alpha(a_s)$. ╥юуфр $\alpha(m)$  ты хЄё  тхъЄюЁюь яюышышэхщэющ ёЄхяхэш ьюэюьр $m$ юЄэюёшЄхы№эю яхЁхьхээ√ї $x_1,\ldots,x_n$ ш фшЇЇхЁхэЎшры№э√ї юяхЁрЄюЁют  $\delta_1,\ldots,\delta_m$. ╤єььє ъюьяюэхэЄ тхъЄюЁр $\alpha(m)$ эрчютхь \textit{ёЄхяхэ№■} ьюэюьр $m$ ш юсючэрўшь ўхЁхч $\deg(m)$.

┴юыхх Єюую, фы  ы■сюую $w\in \mathbb{Z}^{n+m}$ ьюцэю юяЁхфхышЄ№ $w$-ёЄхяхээє■ ЇєэъЎш  $\deg_w$ ъръ $\deg_w(m)=w\cdot \alpha(m)$, уфх $\cdot$ ючэрўрхЄ юс√ўэюх ёъры Ёэюх яЁюшчтхфхэшх. ▀ёэю, ўЄю $\deg_w$ ёютярфрхЄ ё $\deg$ хёыш тёх ъюьяюэхэЄ√ $w$ Ёртэ√ $1$. ┼ёыш яхЁт√х $n$ ъюьяюэхэЄ√  $w$ Ёртэ√ $1$ ш юёЄры№э√х Ёртэ√ $0$, Єю $\deg_w$  ты хЄё  юс∙хщ ёЄхяхэ№■ яю яхЁхьхээ√ь  $x_1,\ldots,x_n$. ╥ръшь юсЁрчюь, ы■сюх $w\in \mathbb{Z}^{n+m}$ юяЁхфхы хЄ уЁрфєшЁютъє
\bes
C=\bigoplus_{i\in \mathbb{Z}} C_i
\ees
рыухсЁ√ $C=R\{x_1,x_2,\ldots,x_n\}$, уфх $C_i$  ты хЄё  $R$-юсюыюўъющ ьюэюьют $w$-ёЄхяхэш  $i$. ╩рцф√щ эхэєыхтющ ¤ыхьхэЄ $c\in C$ юфэючэрўэю яЁхфёЄрты хЄё  т тшфх
\bes
c=c_{i_1}+c_{i_2}+\ldots+c_{i_s}, \ \ i_1<i_2<\ldots<i_s, \ \  0\neq c_{i_j}\in C_{i_j}.
\ees
▌ыхьхэЄ $c_{i_s}$ эрч√трхЄё  {\em ёЄрЁ°хщ юфэюЁюфэющ ўрёЄ№■} ¤ыхьхэЄр $c$ яю юЄэю°хэш■ ъ $w$-ёЄхяхэш $\deg_w$. ╫хЁхч $\overline{c}$ сєфхь юсючэрўрЄ№ ёЄрЁ°є■ юфэюЁюфэє■ ўрёЄ№ $c$ яю юЄэю°хэш■ ъ ЇєэъЎшш ёЄхяхэш $\deg$.

╧єёЄ№ $k$ -- яЁюшчтюы№эюх фшЇЇхЁхэЎшры№эюх яюых їрЁръЄхЁшёЄшъш $0$ ш $B=k\{X\}=k\{x_1,\ldots,x_n \}$ -- рыухсЁр фшЇЇхЁхэЎшры№э√ї ьэюуюўыхэют эрф яюыхь $k$ юЄ ьэюцхёЄтр яхЁхьхээ√ї $X$. ─ы  ы■с√ї $0\neq f, g \in B$, шьххь
$$
\a(fg)=\a(f)+\a(g), \ \
\deg(fg)=\deg(f)+\deg(g), \ \
\overline{fg}=\overline{f} \overline{g}.
$$

▌ыхьхэЄ $f\in B$ эрч√трхЄё  \textit{фшЇЇхЁхэЎшры№эю-рыухсЁршўхёъшь} эрф $k$, хёыш эрщфхЄё  эхэєыхтющ ¤ыхьхэЄ $g\in k\{z\}$ Єръющ, ўЄю $g(f)=0$. ╚эрўх $f\in B$ эрч√трхЄё  \textit{фшЇЇхЁхэЎшры№эю-ЄЁрэёЎхэфхэЄэ√ь} эрф $k$. ▌ыхьхэЄ√ $f_1,f_2,\ldots,f_s\in B$ эрч√тр■Єё  \textit{фшЇЇхЁхэЎшры№эю-рыухсЁршўхёъш чртшёшь√ьш} эрф $k$, хёыш эрщфхЄё  эхэєыхтющ ¤ыхьхэЄ $g\in k\{z_1,\ldots,z_s\}$ Єръющ, ўЄю
$g(f_1,f_2,\ldots,f_s)=0$. ┼ёыш $f_1,f_2,\ldots,f_s$ -- фшЇЇхЁхэЎшры№эю-рыухсЁршўхёъш эхчртшёшь√, Єю уюьюьюЁЇшчь $k\{z_1,\ldots,z_s\} \rightarrow k\{f_1,\ldots,f_s\}$ юяЁхфхыхээ√щ яЁртшыюь $z_i \mapsto f_i$  ты хЄё  шчюьюЁЇшчьюь.

\begin{lemma}\label{l1} ╦■сющ ¤ыхьхэЄ рыухсЁ√ $B=k\{x_1,\ldots,x_n\}$ эх яЁшэрфыхцр∙шщ яюы■ $k$  ты хЄё  фшЇЇхЁхэЎшры№эю-ЄЁрэёЎхэфхэЄэ√ь эрф $k$.
\end{lemma}

\Proof ╙ЄтхЁцфхэшх ыхьь√  ты хЄё  эхёыюцэ√ь ёыхфёЄтшхь шчтхёЄэ√ї ЄхюЁхь ю фшЇЇхЁхэЎшры№эющ ёЄхяхэш ЄЁрэёЎхэфхэЄэюёЄш \cite[├ыртр 2]{Kol}. ╠√ чфхё№ яЁхфырурхь яЁ ьюх фюърчрЄхы№ёЄтю, шёяюы№чє  юс√ўэє■ рыухсЁршўхёъє■ чртшёшьюёЄ№ ¤ыхьхэЄют.

─ы  ы■с√ї $u,v\in X^{\Theta}$ яюыюцшь $u < v$, хёыш $\deg(u)<\deg(v)$ шыш $\deg(u)=\deg(v)$ ш
$\a(u) <\a(v)$ юЄэюёшЄхы№эю ыхъёшъюуЁрЇшўхёъюую яюЁ фър т $\mathbb{Z}_+^{n+m}$.

╧єёЄ№ $0\neq f\in B$. ╧єёЄ№ $u$ -- эршсюы№°шщ ¤ыхьхэЄ шч $X^{\Theta}$, ъюЄюЁ√щ яЁшёєЄёЄтєхЄ т чряшёш $f$. ╥ръющ ¤ыхьхэЄ $u$ эрч√трхЄё  {\em ышфхЁюь} $f$ юЄэюёшЄхы№эю чрфрээюую яюЁ фър $\leq$ эр $X^{\Theta}$ \cite[├ыртр 1]{Kol}. ╦хуъю яюэ Є№, ўЄю ышфхЁюь ¤ыхьхэЄр $f^{\theta}$  ты хЄё  $u^{\theta}$, Є.х. $u^{\Theta}$  ты хЄё  ьэюцхёЄтюь ышфхЁют ьэюцхёЄтр ¤ыхьхэЄют $f^{\Theta}$.

╧юыюцшь $W=X^{\Theta}\setminus u^{\Theta}$. ╥юуфр ьэюцхёЄтю тёхї ¤ыхьхэЄют $u^{\Theta}$  ты хЄё  рыухсЁршўхёъш эхчртшёшь√ь эрф $k[W]$, Єръ ъръ $u^{\Theta}$ ш $W$ юяЁхфхы хЄ Ёрчсшхэшх ьэюцхёЄтр $X^{\Theta}$, ъюЄюЁюх рыухсЁршўхёъш эхчртшёшью эрф $k$.

╬ЄьхЄшь, ўЄю $f$  ты хЄё  фшЇЇхЁхэЎшры№эю-рыухсЁршўхёъшь эрф $k$ Єюуфр ш Єюы№ъю Єюуфр, ъюуфр ьэюцхёЄтю ¤ыхьхэЄют $f^{\Theta}$ -- рыухсЁршўхёъш чртшёшью эрф $k$. ╦■ср  рыухсЁршўхёър  чртшёшьюёЄ№ ¤ыхьхэЄют $f^{\Theta}$ эрф $k$ тхфхЄ ъ рыухсЁршўхёъющ чртшёшьюёЄш $u^{\Theta}$ эрф $k[W]$, ўЄю эхтючьюцэю. $\Box$
\\

┼ёыш $f_1,f_2,\ldots,f_r\in B$, Єю ўхЁхч $k\{f_1,f_2,\ldots,f_r\}$ сєфхь юсючэрўрЄ№ яюфрыухсЁє $B$ яюЁюцфхээє■ ¤ыхьхэЄрьш $f_1,f_2,\ldots,f_r$. ╬ЄьхЄшь, ўЄю Єрър  чряшё№ эх ючэрўрхЄ фшЇЇхЁхэЎшры№эє■-рыухсЁршўхёъє■ эхчртшёшьюёЄ№ ¤ыхьхэЄют $f_1,f_2,\ldots,f_r$, Є.х. $k\{f_1,f_2,\ldots,f_r\}$ эх юс чрЄхы№эю шчюьюЁЇэр рыухсЁх фшЇЇхЁхэЎшры№э√ї ьэюуюўыхэют. └эрыюушўэр  чряшё№ ўрёЄю шёяюы№чєхЄё  фы  юсючэрўхэш  яюфрыухсЁ рыухсЁ ьэюуюўыхэют т рЇЇшээющ рыухсЁршўхёъющ ухюьхЄЁшш.
╙ЄтхЁцфхэшх ёыхфє■∙хщ ыхьь√ тхЁэю фы  ы■с√ї юфэюЁюфэ√ї ётюсюфэ√ї рыухсЁ (ёь. эряЁшьхЁ \cite{Shir}).

\begin{lemma}\label{l2} ╧єёЄ№ $f_1,f_2,\ldots,f_r\in B$ ш $u\in k\{f_1,f_2,\ldots,f_r\}$. ╥юуфр хёыш $\overline{f_1},\overline{f_2},\ldots,\overline{f_r}$ -- фшЇЇхЁхэЎшры№эю-рыухсЁршўхёъш эхчртшёшь√, Єю $\overline{u}\in k\{\overline{f_1},\overline{f_2},\ldots,\overline{f_r}\}$.
\end{lemma}

\Proof ╧єёЄ№ $u=u(z_1,\ldots,z_r)\in k\{z_1,\ldots,z_r\}$ ш яєёЄ№ Єръцх $\deg (f_i)=n_i$, уфх $1\leq i\leq r$. ╧юыюцшь $w=(n_1,n_2,\ldots,n_r,1,\ldots,1)$ ш ЁрёёьюЄЁшь т рыухсЁх $k\{z_1,\ldots,z_r\}$ ЇєэъЎш■ ёЄхяхэш $\deg_w$. ╥юуфр
$u=u'+\widetilde{u}$, уфх $\widetilde{u}$ -- ёЄрЁ°р  юфэюЁюфэр  ўрёЄ№ $u$ юЄэюёшЄхы№эю $\deg_w$ ш $\deg_w(u')<\deg_w(\widetilde{u})$. ╧єёЄ№ $\deg_w(u)=k$.
╟рьхЄшь, ўЄю $f_i=f_i'+\overline{f_i}$ фы  тёхї $i$.
╥юуфр
$$u(f_1,\ldots,f_r)=u'(f_1,\ldots,f_r)+\widetilde{u}(f_1,\ldots,f_r)
=w'+\widetilde{u}(\overline{f_1},\overline{f_2},\ldots,\overline{f_r}),$$
уфх $\deg(w')<k$. ╥ръ ъръ $\overline{f_1},\overline{f_2},\ldots,\overline{f_r}$ -- фшЇЇхЁхэЎшры№эю-рыухсЁршўхёъш эхчртшёшь√, Єю $\widetilde{u}(\overline{f_1},\overline{f_2},\ldots,\overline{f_r})$ эх Ёртхэ эєы■ ш шьххЄ ёЄхяхэ№ $k$ т ёшыє т√сюЁр $w$. ╤ыхфютрЄхы№эю,
$\overline{u}=\widetilde{u}(\overline{f_1},\overline{f_2},\ldots,\overline{f_r})\in k\{\overline{f_1},\overline{f_2},\ldots,\overline{f_r}\}$.
$\Box$

\begin{corol} \label{ё1} ╧єёЄ№ $0\neq f \in B$. ┼ёыш $a\in k\{f\}$, Єю
$\overline{a}\in k\{\overline{f}\}$.
\end{corol}

\Proof ═хяюёЁхфёЄтхээю т√ЄхърхЄ шч  ыхьь \ref{l1} ш \ref{l2}. $\Box$

\section{└ьры№урьшЁютрээюх ётюсюфэюх яЁюшчтхфхэшх}
\hspace*{\parindent}

╧єёЄ№ $A=k\{x,y\}$ -- рыухсЁр фшЇЇхЁхэЎшры№э√ї ьэюуюўыхэют юЄ фтєї яхЁхьхээ√ї  $x,y$ ш яєёЄ№ $Aut(A)$ -- уЁєяяр ртЄюьюЁЇшчьют рыухсЁ√ $A$.

╫хЁхч $\varphi=(f_1,f_2)$ юсючэрўшь ртЄюьюЁЇшчь рыухсЁ√ $A$ Єръющ, ўЄю $\varphi(x)=f_1, \varphi(y)=f_2$.
└тЄюьюЁЇшчь√ тшфр
$$\sigma(1,a,f)=(a x+f(y),y),$$
$$\sigma(2,a,g)=(x,a y+g(x)),$$
уфх $0\neq a \in k$, $f(y)\in k\{y\}$, $g(x)\in k\{x\}$, эрч√тр■Єё  \textit{¤ыхьхэЄрЁэ√ьш}. ╧юфуЁєяяр $T(A)$ уЁєяя√ $Aut(A)$, яюЁюцфхээр  тёхьш ¤ыхьхэЄрЁэ√ьш ртЄюьюЁЇшчьрьш, эрч√трхЄё  \textit{яюфуЁєяяющ Ёєўэ√ї ртЄюьюЁЇшчьют}. ═х Ёєўэ√х ртЄюьюЁЇшчь√ эрч√тр■Єё  \textit{фшъшьш}.

─ы  ртЄюьюЁЇшчьр $\theta=(f_1,f_2)\in Aut(A)$ юяЁхфхышь ёЄхяхэ№, яюырур 
$$\deg(\theta)=\deg(f_1)+\deg(f_2).$$

┼ёыш
$$\theta=(f_1,f_2),\:\: \varphi=(g_1,g_2),$$
Єю яЁюшчтхфхэшх т $Aut(A)$ юяЁхфхы хЄё  ёыхфє■∙хщ ЇюЁьєыющ:
$$\theta\circ \varphi=(g_1(f_1,f_2),g_2(f_1,f_2)).$$

╧єёЄ№ $Af_2(A)$ -- уЁєяяр рЇЇшээ√ї ртЄюьюЁЇшчьют рыухсЁ√ $A$, Є.х. уЁєяяр ртЄюьюЁЇшчьют тшфр $(a_1x+b_1y+c_1,a_2x+b_2y+c_2)$, уфх $a_i,b_i,c_i\in k, a_1b_2\neq a_2b_1$, $Tr_2(A)$ -- уЁєяяр ЄЁхєуюы№э√ї ртЄюьюЁЇшчьют рыухсЁ√ $A$, Є.х. уЁєяяр ртЄюьюЁЇшчьют тшфр $(ax+f(y),by+c)$, уфх $0\neq a,b\in k,\: c\in k, \:f(y)\in k\{y\}$, ш яєёЄ№ $C=Af_2(A)\cap Tr_2(A)$.

╧єёЄ№ $G$ -- яЁюшчтюы№эр  уЁєяяр, $G_0, G_1, G_2$ -- яюфуЁєяя√ уЁєяя√ $G$, яЁшўхь $G_0=G_1\cap G_2$. ├Ёєяяр $G$ эрч√трхЄё  \textit{ётюсюфэ√ь яЁюшчтхфхэшхь яюфуЁєяя $G_1$ ш $G_2$ ё юс·хфшэхээющ яюфуЁєяяющ $G_0$} ш юсючэрўрхЄё  $G=G_1\ast_{G_0} G_2$, хёыш
\begin{enumerate}
\item[(a)] $G$ яюЁюцфрхЄё  яюфуЁєяярьш $G_1$ ш $G_2$;
\item[(b)] ╬яЁхфхы ■∙шх ёююЄэю°хэш  уЁєяя√ $G$ ёюёЄю Є Єюы№ъю шч юяЁхфхы ■∙шї ёююЄэю°хэшщ яюфуЁєяя $G_1$ ш $G_2$.
\end{enumerate}

┼ёыш $S_1$ -- ёшёЄхьр ыхт√ї яЁхфёЄртшЄхыхщ $G_1$ яю $G_0$, $S_2$ -- ёшёЄхьр ыхт√ї яЁхфёЄртшЄхыхщ $G_2$ яю $G_0$, Єю уЁєяяр $G$  ты хЄё  ётюсюфэ√ь яЁюшчтхфхэшхь яюфуЁєяя $G_1$ ш $G_2$ ё юс·хфшэхээющ яюфуЁєяяющ $G_0$ (ёь. эряЁшьхЁ \cite{Mag}) т Єюь ш Єюы№ъю т Єюь ёыєўрх, ъюуфр ърцф√щ $g\in G$ юфэючэрўэю яЁхфёЄрты хЄё  т тшфх
$$g=g_1\ldots g_k c,$$
уфх $g_i\in S_1 \cup S_2,\: i=1,\ldots,k,$ $g_i, g_{i+1}$ юфэютЁхьхээю эх яЁшэрфыхцрЄ $S_1$ шыш $S_2$, $c\in G_0$.

╟ряшё№ $h_i(y)$ т фюърчрЄхы№ёЄтрї ёыхфє■∙шї эхёъюы№ъшї ыхьь ючэрўрхЄ, ўЄю $h_i(y)\in k\{y\}$ -- юфэюЁюфэ√щ фшЇЇхЁхэЎшры№э√щ ьэюуюўыхэ ёЄхяхэш $i$ яю юЄэю°хэш■ ъ ЇєэъЎшш ёЄхяхэш $\deg$ юЄ юфэющ яхЁхьхээющ $y$. ▀ёэю, ўЄю
 $h_0(y)\in k$.

\begin{lemma} \label{l3} a) ╤шёЄхьр ¤ыхьхэЄют
$$A_0=\{id=(x,y), \gamma=(y,x+ay) | a\in k\}$$
 ты хЄё  ёшёЄхьющ яЁхфёЄртшЄхыхщ ыхт√ї ёьхцэ√ї ъырёёют $Af_2(A)$ яю яюфуЁєяях $C$.

b) ╤шёЄхьр ¤ыхьхэЄют
$$B_0=\{\beta=(x+q(y),y) | q(y)=h_n(y)+\ldots+h_2(y) \}$$
 ты хЄё  ёшёЄхьющ яЁхфёЄртшЄхыхщ ыхт√ї ёьхцэ√ї ъырёёют $Tr_2(A)$ яю яюфуЁєяях $C$.
\end{lemma}

\Proof ╧ЁютхЁшь єёыютшх a). ╧єёЄ№ $l\in Af_2(A)$. ╠√ фюыцэ√ яюърчрЄ№, ўЄю фы  ы■сюую $l$ эрщфєЄё  $\gamma \in A_0,\:\eta\in C$ Єръшх, ўЄю $l=\gamma \circ \eta$.

┼ёыш $l=(a_1 x+b_1 y+c_1, a_2 x+b_2 y+ c_2)$, уфх $a_2\neq 0$, Єю яюыюцшь  $\gamma=(y,x+\frac{b_2}{a_2} y)$, $\eta=((b_1-\frac{a_1 b_2}{a_2})x+a_1 y+c_1, a_2 y+c_2)$. ╥юуфр $l$ яЁхфёЄрты хЄё  т тшфх
$$l=(y,x+\frac{b_2}{a_2} y)\circ ((b_1-\frac{a_1 b_2}{a_2})x+a_1 y+c_1, a_2 y+c_2)=\gamma \circ \eta.$$

┼ёыш $a_2=0$, Єю $\gamma=id$, $\eta=l$, Є.х. $l=id\circ l$.

─юяєёЄшь $\gamma_1=(y,x+a_1 y)$, $\gamma_2=(y, x+a_2 y)$ ш $\gamma_1 C=\gamma_2 C$, Єюуфр
$$\gamma_1^{-1}\circ \gamma_2=(-a_1 x+y, x)\circ (y, x+a_2y)=(x, (-a_1+a_2)x+y).$$
╬Єё■фр ёыхфєхЄ, ўЄю $\gamma_1^{-1}\circ \gamma_2\in C$ Єюуфр ш Єюы№ъю Єюуфр, ъюуфр $a_1=a_2$. ╤ыхфютрЄхы№эю, $\gamma_1=\gamma_2$.

╥хяхЁ№ яЁютхЁшь єёыютшх b). ╧єёЄ№ $\psi=(ax+h(y),by+c)\in Tr_2(A)$ ш  $h(y)=h_n(y)+\ldots+h_1(y)+h_0(y)$.
╠√ фюыцэ√ яюърчрЄ№, ўЄю фы  ы■сюую $\psi$ эрщфєЄё  $\beta\in B_0,\:\mu\in C$ Єръшх, ўЄю $\psi=\beta\circ \mu$. ╧юыюцшь $\beta=(x+q(y),y)$, $\mu=(ax+h_1(y)+h_0(y),by+c)$, уфх $q(y)=h_n(y)+\ldots+h_2(y)$. ╥юуфр $\psi$ яЁхфёЄрты хЄё  т тшфх
$$\psi=(x+\frac{1}{a} q(y),y) \circ(ax+h_1(y)+h_0(y),by+c)=\beta\circ \mu.$$

─юяєёЄшь, $\beta_1=(x+q(y),y),$ $\beta_2=(x+q^{(1)}(y),y)$ ш $\beta_1C=\beta_2C$. ╥юуфр шьххь
$$\beta_1^{-1}\circ\beta_2=(x-q(y),y)\circ (x+q^{(1)}(y),y)=(x-q(y)+q^{(1)}(y),y).$$
╬Єё■фр ёыхфєхЄ, ўЄю $\beta_1^{-1}\circ \beta_2\in C$ Єюуфр ш Єюы№ъю Єюуфр, ъюуфр $q(y)=q^{(1)}(y)$. ╤ыхфютрЄхы№эю, $\beta_1=\beta_2$.
$\Box$

\begin{lemma} \label{l4} ╧єёЄ№ $A_0, B_0$ -- ьэюцхёЄтр, юяЁхфхыхээ√х т ыхььх \ref{l3}. ╥юуфр ы■сющ Ёєўэющ ртЄюьюЁЇшчь $\varphi$ рыухсЁ√ $A$ ЁрчырурхЄё  т яЁюшчтхфхэшх тшфр
\begin{gather} \label{varphi}
\varphi=\gamma_1\circ \beta_1 \circ \gamma_2 \circ \beta_2 \circ\ldots \circ \gamma_k \circ\beta_k\circ \gamma_{k+1}\circ \lambda,
\end{gather}
уфх $\gamma_i\in A_0,\:\gamma_2,\ldots,\gamma_k\neq id$, $\beta_i\in B_0,\:\beta_1,\ldots,\beta_k\neq id$, $\lambda\in C$.
\end{lemma}

\Proof ╬ўхтшфэю, ўЄю
$$(ax+h(y),y)=(x+\frac{1}{a} q(y),y)\circ (ax+h_1(y)+h_0(y),y),$$
уфх $h(y)=h_n(y)+\ldots+h_2(y)+h_1(y)+h_0(y)$, $q(y)=h_n(y)+\ldots+h_2(y)$,
$$(x,by+h^{(1)}(x))=(y,x)\circ (x+\frac{1}{b} q^{(1)}(y),y)\circ (y,bx+h^{(1)}_1(y)+h^{(1)}_0(y)),$$
уфх $h^{(1)}(y)=h^{(1)}_m(y)+\ldots+h^{(1)}_2(y)+h^{(1)}_1(y)+h^{(1)}_0(y)$, $q^{(1)}(y)=h^{(1)}_m(y)+\ldots+h^{(1)}_2(y)$,
Є.х. ы■сющ ¤ыхьхэЄрЁэ√щ ртЄюьюЁЇшчь шьххЄ тшф
$$l_1\circ \beta \circ l_2,$$
уфх $\beta\in B_0$, $l_1, l_2\in Af_2(A)$.

╦■сющ Ёєўэющ ртЄюьюЁЇшчь $\varphi$ яЁхфёЄрты хЄё  т тшфх ъюьяючшЎшш ¤ыхьхэЄрЁэ√ї ртЄюьюЁЇшчьют $\varphi_1,\varphi_2,\ldots,\varphi_n$, Є.х.
$$\varphi=\varphi_1 \circ \varphi_2 \circ \ldots \circ \varphi_n.$$
╤ыхфютрЄхы№эю, шьххь
\begin{gather} \label{razlag}
\varphi=l_1 \circ \beta_1 \circ l_2 \circ \beta_2 \circ \ldots \circ l_n \circ \beta_n \circ l_{n+1},
\end{gather}
уфх $\beta_i\in B_0$, $l_i\in Af_2(A)$.

─юърцхь шэфєъЎшхщ яю $n$, ўЄю $\varphi$ яЁхфёЄрты хЄё  т тшфх яЁюшчтхфхэш  \eqref{varphi}, ё $k\leq n$.

╤юуырёэю ыхььх \ref{l3} ртЄюьюЁЇшчь $l_1$ чряшё√трхЄё  т тшфх $\gamma_1\circ \lambda_1$, уфх $\gamma_1\in A_0$, $\lambda_1\in C$.
╥юуфр
$$l_1\circ\beta_1=\gamma_1\circ\lambda_1\circ\beta_1.$$

╧єёЄ№ $\lambda_1=(ax+by+c,b_1y+c_1)$, $\beta_1=(x+q(y),y)$. ╥юуфр
$$\lambda_1\circ \beta_1\circ \lambda_1^{-1}=(x+\frac{1}{a}q(b_1y+c_1),y).$$
╫хЁхч $q_{<2}(b_1y+c_1)$ юсючэрўшь ышэхщэє■ ўрёЄ№ фшЇЇхЁхэЎшры№эюую ьэюуюўыхэр $q(b_1y+c_1)$. ╧єёЄ№ $\lambda=(x-\frac{1}{a}q_{<2}(b_1y+c_1),y)$. ▀ёэю, ўЄю $\lambda\in C$ ш $\lambda_1^{-1}\circ\lambda\in C$. ╬сючэрўшь $\lambda_1^{-1}\circ\lambda$ ўхЁхч $\lambda_2^{-1}$. ╥юуфр
$$l_1\circ\beta_1=\gamma_1\circ\lambda_1\circ\beta_1=\gamma_1\circ\beta'_1\circ\lambda_2,$$
уфх $\beta'_1=\lambda_1\circ\beta_1\circ\lambda_2^{-1}=(x+\frac{1}{a}q(b_1y+c_1)-\frac{1}{a}q_{<2}(b_1y+c_1),y)\in B_0$.
╚ьххь
$$\varphi=\gamma_1 \circ \beta'_1 \circ (\lambda_2 \circ l_2) \circ \beta_2 \circ \ldots \circ l_n \circ \beta_n \circ l_{n+1}.$$
╧ю шэфєъЄштэюьє яЁхфяюыюцхэш■, яЁюшчтхфхэшх
$$(\lambda_2 \circ l_2)\circ \beta_2 \circ \ldots \circ l_n\circ \beta_n\circ l_{n+1}$$
чряшё√трхЄё  т тшфх $$\gamma_2\circ\beta'_2\circ\gamma_3\circ\ldots\circ\gamma_k\circ\beta'_k\circ\gamma_{k+1}\circ\lambda,\:k\leq n.$$
╤ыхфютрЄхы№эю,
$$\varphi=\gamma_1\circ\beta'_1\circ\gamma_2\circ\beta'_2\circ\ldots\circ\gamma_k\circ\beta'_k\circ\gamma_{k+1}\circ\lambda.$$
┼ёыш $\gamma_2\neq id$, Єю яюыєўхээюх яЁхфёЄртыхэшх шьххЄ тшф \eqref{varphi}. ╥хяхЁ№ ЁрёёьюЄЁшь ёыєўрщ ъюуфр $\gamma_2=id$. ╥ръ ъръ $\beta'_1\circ\beta'_2=\beta''_2\in B_0$, Єю
$$\varphi=\gamma_1\circ\beta'_1\circ\beta'_2\circ\gamma_3\circ\ldots\circ\gamma_k\circ\beta'_k\circ\gamma_{k+1}\circ\lambda=\gamma_1\circ\beta''_2\circ\gamma_3\circ\ldots\circ\gamma_k\circ\beta'_k\circ\gamma_{k+1}\circ\lambda.$$
╧юёъюы№ъє $k-1<n$, Єю яю шэфєъЄштэюьє яЁхфяюыюцхэш■ $\varphi$ чряшё√трхЄё  т тшфх \eqref{varphi}.
$\Box$

\begin{lemma} \label{l5} ╧єёЄ№ $\varphi=(f_1,f_2)$ -- ртЄюьюЁЇшчь рыухсЁ√ $A$, яЁхфёЄртшь√щ т тшфх яЁюшчтхфхэш 
$$\varphi=(f_1,f_2)=\beta_1\circ\gamma_2\circ\beta_2\circ\ldots\circ\gamma_k\circ\beta_k,$$
уфх $id\neq \gamma_i\in A_0$, $id\neq \beta_i\in B_0$ фы  тёхї $i$.
┼ёыш $\beta_i=(x+q_i(y),y)$, $\deg(q_i(y))=n_i$ ш $s_i$ -- ёЄхяхэ№ $q_i(y)$ яю яхЁхьхээющ $y$ фы  тёхї $1\leq i \leq k$, Єю
$$\deg(f_1)=n_k+(n_{k-1}-1)s_k+\ldots+(n_1-1)s_k s_{k-1}\ldots s_2,$$
$$\deg(f_2)=n_{k-1}+(n_{k-2}-1)s_{k-1}+\ldots+(n_1-1)s_{k-1} s_{k-2}\ldots s_2, \:\:\text{хёыш}\:\: k>1$$
ш
$$\deg(f_2)=1,\:\: \text{хёыш}\:\: k=1.$$
\end{lemma}

\Proof ╙ЄтхЁцфхэшх ыхьь√ фюърцхь шэфєъЎшхщ яю $k$. ┼ёыш $k=1$, Єю $\varphi=\beta_1$ ш
$$\deg(f_1)=\deg(q_1(y))=n_1,$$
$$\deg(f_2)=1.$$
╧Ёхфяюыюцшь, ўЄю єЄтхЁцфхэшх ыхьь√ т√яюыэ хЄё  фы  $k-1$. ╧юыюцшь,
$$\varphi_1=\beta_1\circ\gamma_2\circ\beta_2\circ\ldots\circ\gamma_{k-1}\circ\beta_{k-1}=(g_1,g_2).$$
╧ю шэфєъЄштэюьє яЁхфяюыюцхэш■, шьххь
$$\deg(g_1)=n_{k-1}+(n_{k-2}-1)s_{k-1}+\ldots+(n_1-1)s_{k-1} s_{k-2}\ldots s_2,$$
$$\deg(g_2)=n_{k-2}+(n_{k-3}-1)s_{k-2}+\ldots+(n_1-1)s_{k-2} s_{k-3}\ldots s_2.$$
╥юуфр
$$\varphi=(f_1,f_2)=\beta_1\circ\gamma_2\circ\beta_2\circ\ldots\circ\gamma_k\circ\beta_k=\varphi_1\circ\gamma_k\circ\beta_k=(g_1,g_2)\circ\gamma_k\circ\beta_k.$$
╧Ёшьхэ   $\gamma_k=(y,x+ay)$ ъ $(g_1,g_2)$, яюыєўшь
$$(u_1,u_2)=(g_1,g_2)\circ\gamma_k=(g_2,g_1+a g_2).$$
╥юуфр
$$\deg(u_1)=\deg(g_2)=n_{k-2}+(n_{k-3}-1)s_{k-2}+\ldots+(n_1-1)s_{k-2} s_{k-3}\ldots s_2,$$
\bes
\deg(u_2)=max\{\deg(g_1),\deg(g_2)\}\\
=n_{k-1}+(n_{k-2}-1)s_{k-1}+\ldots+(n_1-1)s_{k-1} s_{k-2}\ldots s_2.
\ees
─рыхх,
$$\varphi=(f_1,f_2)=(u_1,u_2)\circ\beta_k=(u_1,u_2)\circ(x+q_k(y),y)=(u_1+q_k(u_2),u_2).$$
╤ыхфютрЄхы№эю,
$$\deg(f_1)=max\{\deg(u_1),\deg(q_k(u_2))\},$$
$$\deg(f_2)=\deg(u_2).$$
═ряюьэшь, ўЄю $\deg(q_k)=n_k$ ш
$$\deg(u_2)=n_{k-1}+(n_{k-2}-1)s_{k-1}+\ldots+(n_1-1)s_{k-1} s_{k-2}\ldots s_2.$$
 ╟рьхЄшь, ўЄю
$$\overline{q_k(u_2)}=\widetilde{q_k}(\overline{u_2}),$$
уфх $\widetilde{q_k}$ -- ёЄрЁ°р  юфэюЁюфэр  ўрёЄ№ $q_k$ юЄэюёшЄхы№эю $\deg_w$, $w=(t,\underbrace{1,1,\ldots,1}_{m})$ ш $t=\deg(u_2)$. ╥юуфр
$$\deg(q_k(u_2))=\deg(\overline{q_k(u_2)})=\deg(\widetilde{q_k}(\overline{u_2}))=\deg_w(q_k)=(t,1,1,\ldots,1)\cdot \alpha(q_k)$$
$$=\deg(q_k)+(t-1)s_k=n_k+(n_{k-1}-1)s_k+(n_{k-2}-1)s_k s_{k-1}+\ldots+(n_1-1)s_k s_{k-1}\ldots s_2.$$
╤ыхфютрЄхы№эю,
$$\deg(f_1)=n_k+(n_{k-1}-1)s_k+\ldots+(n_1-1)s_k s_{k-1}\ldots s_2,$$
$$\deg(f_2)=n_{k-1}+(n_{k-2}-1)s_{k-1}+\ldots+(n_1-1)s_{k-1} s_{k-2}\ldots s_2.$$
╫Єю ш ЄЁхсютрыюё№ фюърчрЄ№.
$\Box$

\begin{lemma} \label{l6} ╨рчыюцхэшх \eqref{varphi} ртЄюьюЁЇшчьр $\varphi$ шч ыхьь√ \ref{l4}  ты хЄё  юфэючэрўэ√ь.
\end{lemma}

\Proof ─юёЄрЄюўэю яюърчрЄ№, ўЄю
$$\gamma_1\circ\beta_1\circ\gamma_2\circ\beta_2\circ\ldots\circ\gamma_k\circ \beta_k\circ\gamma_{k+1}\circ\lambda\neq id,$$
яЁш $k\geq 1$, $\gamma_i\in A_0,\:\gamma_2,\ldots,\gamma_k\neq id$, $\beta_i\in B_0,\:\beta_1,\ldots,\beta_k\neq id$, $\lambda\in C$.

─юърцхь юЄ яЁюЄштэюую. ─юяєёЄшь
$$\gamma_1\circ\beta_1\circ\gamma_2\circ\beta_2\circ\ldots\circ\gamma_k\circ\beta_k\circ\gamma_{k+1}\circ\lambda= id.$$
╥юуфр
\begin{gather} \label{prot}
\beta_1\circ\gamma_2\circ\beta_2\circ\ldots\circ\gamma_k\circ\beta_k=\gamma_1^{-1}\circ\lambda^{-1}\circ\gamma_{k+1}^{-1}.
\end{gather}
╤юуырёэю ыхььх \ref{l5} ртЄюьюЁЇшчь
$$\varphi=(f_1,f_2)=\beta_1\circ\gamma_2\circ\beta_2\circ\ldots\circ\gamma_k\circ\beta_k$$
шьххЄ ёЄхяхэ№
$$\deg(\varphi)=\deg(f_1)+\deg(f_2)=n_k+(n_{k-1}-1)s_k+\ldots+(n_1-1)s_k s_{k-1}\ldots s_2$$
$$+n_{k-1}+(n_{k-2}-1)s_{k-1}+\ldots+(n_1-1)s_{k-1} s_{k-2}\ldots s_2.$$
╧Ёртє■ ўрёЄ№ ЁртхэёЄтр \eqref{prot} юсючэрўшь ўхЁхч $\rho$, Є.х.
$$\rho=\gamma_1^{-1}\circ\lambda^{-1}\circ\gamma_{k+1}^{-1}.$$
▀ёэю, ўЄю $\rho \in Af_2(A)$ ш $\deg(\rho)=2$.
╤ыхфютрЄхы№эю, $\deg(\varphi)\neq \deg(\rho)$, ўЄю яЁюЄштюЁхўшЄ ЁртхэёЄтє \eqref{prot}.
$\Box$

\begin{theorem} \label{t1} ├Ёєяяр Ёєўэ√ї ртЄюьюЁЇшчьют рыухсЁ√ $A=k\{x,y\}$  ты хЄё  ётюсюфэ√ь яЁюшчтхфхэшхь яюфуЁєяя рЇЇшээ√ї ртЄюьюЁЇшчьют $Af_2(A)$ ш ЄЁхєуюы№э√ї ртЄюьюЁЇшчьют $Tr_2(A)$ ё юс·хфшэхээющ яюфуЁєяяющ $C=Af_2(A)\cap Tr_2(A)$, Є.х.
$$T(A)=Af_2(A) \ast_{C} Tr_2(A).$$
\end{theorem}

\Proof ╥ръ ъръ $A_0$ ш $B_0$ -- ёшёЄхь√ ыхт√ї ёьхцэ√ї ъырёёют $Af_2(A)$ ш $Tr_2(A)$ яю яюфуЁєяях $C$, Єю яю ыхььх \ref{l4} ш яю ыхььх \ref{l6} ы■сющ Ёєўэющ ртЄюьюЁЇшчь юфэючэрўэю яЁхфёЄрты хЄё  т тшфх \eqref{varphi}. ╤юуырёэю \cite{Mag},
$$T(A)=Af_2(A) \ast_{C} Tr_2(A). \ \ \Box $$

\section{╤юъЁрЄшьюёЄ№ Ёєўэ√ї ртЄюьюЁЇшчьют}
\hspace*{\parindent}

═ряюьэшь, ўЄю $\overline{f}$ -- ёЄрЁ°р  юфэюЁюфэр  ўрёЄ№ $f$ яю юЄэю°хэш■ ъ ЇєэъЎшш ёЄхяхэш $\deg$ ш ёЄхяхэ№ ртЄюьюЁЇшчьр $\theta=(f_1,f_2)$ юяЁхфхы хЄё  ёыхфє■∙шь ёяюёюсюь:
$$\deg(\theta)=\deg(f_1)+\deg(f_2).$$

╧ЁхюсЁрчютрэшх $(f_1,f_2)$, ъюЄюЁюх чрьхэ хЄ Єюы№ъю юфшэ ¤ыхьхэЄ $f_i\:(i=1,2)$ эр ¤ыхьхэЄ тшфр $\alpha f_i+g$, уфх $0 \neq \alpha \in k$, $g\in k\{f_j |j \neq i\}$, эрч√трхЄё  \textit{¤ыхьхэЄрЁэ√ь}.

╟ряшё№ $\theta \rightarrow \varphi$ ючэрўрхЄ, ўЄю $\varphi$ яюыєўрхЄё  шч $\theta$ ё яюью∙№■ юфэюую ¤ыхьхэЄрЁэюую яЁхюсЁрчютрэш . └тЄюьюЁЇшчь $\theta$ эрч√трхЄё  \textit{¤ыхьхэЄрЁэю ёюъЁрЄшь√ь}, хёыш ёє∙хёЄтєхЄ ртЄюьюЁЇшчь $\varphi$ Єръющ, ўЄю $\theta \rightarrow \varphi$ ш $\deg(\varphi)<\deg(\theta)$.

\begin{lemma} \label{l7} ╧єёЄ№ $\theta=(f_1,f_2)$ -- эх рЇЇшээ√щ Ёєўэющ ртЄюьюЁЇшчь рыухсЁ√ $A=k\{x,y\}$. ┼ёыш $\overline{f_1}$ ш $\overline{f_2}$ -- ышэхщэю чртшёшь√, Єю ртЄюьюЁЇшчь $\pi$  ты хЄё  ¤ыхьхэЄрЁэю ёюъЁрЄшь√ь.
\end{lemma}

\Proof ╧єёЄ№ $\overline{f_1}= \gamma\overline{f_2}$.
╨рёёьюЄЁшь ¤ыхьхэЄрЁэюх яЁхюсЁрчютрэшх $$\theta=(f_1,f_2)\rightarrow (f_1-\gamma f_2,f_2)=\sigma,$$
уфх $\gamma\in k^\ast$. ╚ьххь $\deg(f_1)>\deg(f_1-\gamma f_2)$. ╬Єё■фр ёыхфєхЄ, ўЄю $\deg(\theta)>\deg(\sigma)$ ш ртЄюьюЁЇшчь $\pi$  ты хЄё  ¤ыхьхэЄрЁэю ёюъЁрЄшь√ь.
$\Box$

\begin{theorem} \label{t2} ╦■сющ эх рЇЇшээ√щ Ёєўэющ ртЄюьюЁЇшчь рыухсЁ√ $A=k\{x,y\}$  ты хЄё  ¤ыхьхэЄрЁэю ёюъЁрЄшь√ь.
\end{theorem}

\Proof ╧єёЄ№ $\theta=(f_1,f_2)$ -- эх рЇЇшээ√щ Ёєўэющ ртЄюьюЁЇшчь рыухсЁ√ $A$. ╧ю ыхььх \ref{l4} $\theta$ чряшё√трхЄё  т тшфх \eqref{varphi}. ┼ёыш $\gamma_{k+1}\circ \lambda=id$, Єю
$$\theta=\gamma_1\circ \beta_1 \circ \gamma_2 \circ \beta_2 \circ\ldots \circ \gamma_k \circ\beta_k=(f_1,f_2).$$
╧юыюцшь $$\tau=\gamma_1\circ \beta_1 \circ \gamma_2 \circ \beta_2 \circ\ldots \circ \gamma_k=(g_1,g_2).$$
┼ёыш $\beta_k=(x+q_k(y),y)$, Єю
$$\theta=(g_1+q_k(g_2),g_2).$$
╧ю ыхььх \ref{l5}, шьххь
$$\deg(\tau)=\deg(g_1)+\deg(g_2)<\deg(\theta)=\deg(g_1+q_k(g_2))+\deg(g_2).$$
╧юёъюы№ъє
$$\theta\rightarrow \tau,$$
Єю ртЄюьюЁЇшчь $\theta$  ты хЄё  ¤ыхьхэЄрЁэю ёюъЁрЄшь√ь.
─юяєёЄшь, ўЄю
$$\gamma_{k+1}\circ \lambda=(a_1x+b_1y+c_1, a_2x+b_2y+c_2)\neq id.$$
╧юыюцшь
$$\pi=\gamma_1\circ \beta_1 \circ \gamma_2 \circ \beta_2 \circ\ldots \circ \gamma_k \circ\beta_k=(g_1+q_k(g_2),g_2)=(u_1,u_2).$$
╧ю ыхььх \ref{l5} $\deg(u_1)>\deg(u_2)$.

╤ыхфютрЄхы№эю,
$$\theta=\pi\circ \gamma_{k+1}\circ \lambda=(a_1u_1+b_1u_2+c_1,a_2u_1+b_2u_2+c_2)=(f_1,f_2).$$

┼ёыш $a_1,a_2\neq 0$, Єю $\overline{f_1}$ ш $\overline{f_2}$ -- ышэхщэю чртшёшь√ ш яю ыхььх \ref{l7} ртЄюьюЁЇшчь $\theta$  ты хЄё  ¤ыхьхэЄрЁэю ёюъЁрЄшь√ь.

┼ёыш $a_1=0$, Єю $\overline{f_1}=\overline{u_2}$ ш $\overline{f_2}=\overline{u_1}=\overline{q_k(u_2)}$. ┬ ¤Єюь ёыєўрх ртЄюьюЁЇшчь $\theta$ ¤ыхьхэЄрЁэю ёюъЁр∙рхЄё  ё яюью∙№■ ртЄюьюЁЇшчьр $\psi=(f_1,f_2-q_k(f_1))$.

╤ыєўрщ ъюуфр $a_2=0$ рэрыюушўхэ яЁхф√фє∙хьє.
$\Box$

\begin{corol} \label{c2}  ╧єёЄ№ $(f_1,f_2)$ -- эх рЇЇшээ√щ Ёєўэющ ртЄюьюЁЇшчь рыухсЁ√ $A=k\{x,y\}$. ╥юуфр эрщфєЄё  $i$ ш $g\in k\{f_j|j\neq i\}$ Єръшх, ўЄю $\overline{f_i}=\overline{g}$.
\end{corol}

\Proof
┬ ёшыє ЄхюЁхь√ \ref{t2} ртЄюьюЁЇшчь $(f_1,f_2)$  ты хЄё  ¤ыхьхэЄрЁэю ёюъЁрЄшь√ь. ─юяєёЄшь, ўЄю $f_1$  ты хЄё  ёюъЁрЄшь√ь ¤ыхьхэЄюь ¤Єюую ртЄюьюЁЇшчьр. ╥юуфр эрщфхЄё  $g\in k\{f_2\}$ Єръющ, ўЄю $\deg(f_1 -g(f_2))<\deg(f_1)$. ▌Єю ючэрўрхЄ, ўЄю $\overline{f_1}=\overline{g(f_2)}$.
$\Box$

\section{└эрыюу ртЄюьюЁЇшчьр └эшър}
\hspace*{\parindent}

\begin{lemma} \label{l8} ╧єёЄ№ $|\Delta|\geq 2$. ▌эфюьюЁЇшчь $\delta$ рыухсЁ√ $A=k\{x,y\}$ чрфрээ√щ ъръ
$$\delta(x)= x+w^{\delta_2}, \delta(y)= y+w^{\delta_1},$$
уфх $w=x^{\delta_1}-y^{\delta_2}$,  ты хЄё  ртЄюьюЁЇшчьюь.
\end{lemma}

\Proof
╧юыюцшь
$$f_1=x+w^{\delta_2}, f_2=y+w^{\delta_1}.$$
╧юърцхь, ўЄю $k\{x,y\}=k\{f_1,f_2\}$. ╬ўхтшфэю, ўЄю $k\{f_1,f_2\} \subseteq k\{x,y\}$.
╚ьххь
$$x=f_1-w^{\delta_2}, y=f_2-w^{\delta_1}.$$
╤ыхфютрЄхы№эю,
$$w=x^{\delta_1}-y^{\delta_2}={(f_1-w^{\delta_2})}^{\delta_1}-{(f_2-w^{\delta_1})}^{\delta_2}=f_1^{\delta_1}-f_2^{\delta_2}\in k\{f_1, f_2\}$$
ш
$$x=f_1-w^{\delta_2} \in k\{f_1, f_2\}, \:\:y=f_2-w^{\delta_1} \in k\{f_1, f_2\}.$$
▌Єю ючэрўрхЄ, ўЄю $k\{x,y\}\subseteq k\{f_1,f_2\}$.
╬Єё■фр ёыхфєхЄ, ўЄю $\delta$ -- ё■Ё·хъЄштэ√щ уюьюьюЁЇшчь.

╦шэхщэ√х ўрёЄш $f_1$ ш $f_2$ Ёртэ√ $x$ ш $y$, ёююЄтхЄёЄтхээю. ╤ыхфютрЄхы№эю, $f_1$ ш $f_2$ фшЇЇхЁхэЎшры№эю-рыухсЁршўхёъш эхчртшёшь√. ▌Єю яюърч√трхЄ шэ·хъЄштэюёЄ№ уюьюьюЁЇшчьр $\delta$.
$\Box$

\begin{theorem} \label{t3} └тЄюьюЁЇшчь $\delta$ рыухсЁ√ $A=k\{x,y\}$  ты хЄё  фшъшь.
\end{theorem}

\Proof ╚ьххь $f_1, f_2$:
$$\overline{f_1}=\overline{x+x^{\delta_1 \delta_2}-y^{\delta_2^2}}=x^{\delta_1 \delta_2}-y^{\delta_2^2},$$
$$\overline{f_2}=\overline{y+x^{\delta_1^2}-y^{\delta_1 \delta_2}}=x^{\delta_1^2}-y^{\delta_1 \delta_2}.$$

╚ьххь $\deg(x^{\delta_1^2}-y^{\delta_1 \delta_2})=3$ ш $\deg(x^{\delta_1 \delta_2}-y^{\delta_2^2})=3$. ╟рьхЄшь, ўЄю ы■сющ юфэюЁюфэ√щ ¤ыхьхэЄ ёЄхяхэш $3$ рыухсЁ√ $k\{x^{\delta_1 \delta_2}-y^{\delta_2^2}\}$ шьххЄ тшф $a (x^{\delta_1 \delta_2}-y^{\delta_2^2})$ фы  эхъюЄюЁюую $a \in k^*$.
╧ю¤Єюьє $x^{\delta_1^2}-y^{\delta_1 \delta_2} \notin k\{x^{\delta_1 \delta_2}-y^{\delta_2^2}\}$, Єръ ъръ $x^{\delta_1^2}-y^{\delta_1 \delta_2}=a (x^{\delta_1 \delta_2}-y^{\delta_2^2})$ эхтючьюцэю.

└эрыюушўэю, $x^{\delta_1 \delta_2}-y^{\delta_2^2} \notin k\{x^{\delta_1^2}-y^{\delta_1 \delta_2}\}$.

╤ыхфютрЄхы№эю, ртЄюьюЁЇшчь $\delta$ эх єфютыхЄтюЁ хЄ єЄтхЁцфхэш■ ёыхфёЄтш  \ref{c2}, Є.х.  ты хЄё  фшъшь.
$\Box$


\begin{thebibliography}{99}

\bibitem{Jung} Jung H.W.E. Uber ganze birationale Transformationen der Ebene // J. reine angew. Math. -- 1942. -- Vol. 184. -- P. 161--174.

\bibitem{Kulk} Kulk W. Van der. On Polynomial Rings in Two Variables // Nieuw Archief voor Wiskunde. -- 1953. -- Vol. 3, No 1. -- P. 33--41.

\bibitem{Czer} Czerniakiewicz A.G. Automorphisms of a Free Associative Algebra of Rank 2. I, II // Trans. Amer.
Math. Soc. -- 1971. -- Vol. 160. -- P. 393--401; -- 1972. -- Vol. 171. -- P. 309--315.

\bibitem{M-L1} ╠рърЁ-╦шьрэют ╦. └тЄюьюЁЇшчь√ ётюсюфэющ рыухсЁ√ юЄ фтєї яюЁюцфр■∙шї // ╘єэъЎшюэ. рэрышч ш хую яЁшы. -- 1970. -- ╥. 4. -- ╤. 107--108.

\bibitem{M-L2}  Makar-Limanov L., Turusbekova U., Umirbaev U.U. Automorphisms and derivations of free Poisson algebras in two variables // J. Algebra. -- 2009. -- Vol. 322, No 9. -- P. 3318--3330.

\bibitem{Koz}  Kozybaev D.,  Makar-Limanov L., Umirbaev U. The Freiheitssatz and the automorphisms of free right-symmetric algebras // Asian-European Journal of Mathematics. -- 2008. -- Vol. 1. -- P. 243--254.

\bibitem{Cohn64}  Cohn P.M. Subalgebras of free associative algebras // Proc. London Math. Soc. -- 1964. -- Vol. 56. -- P. 618--632.

\bibitem{Lewin}  Lewin J. On Schreier varieties of linear algebras // Trans. Amer. Math. Soc. -- 1968. -- Vol. 132. -- P. 553--562.

\bibitem{14} ╩єЁю° └.├. ═хрёёюЎшрЄштэ√х рыухсЁ√ ш ётюсюфэ√х яЁюшчтхфхэш  рыухсЁ // ╠рЄхь. ёс. -- 1947. -- ╥. 20. -- ╤. 239--262.

\bibitem{15}  ╪шЁ°ют └.╚. ╧юфрыухсЁ√ ётюсюфэ√ї ъюььєЄрЄштэ√ї ш ётюсюфэ√ї рэЄшъюььєЄрЄштэ√ї рыухсЁ // ╠рЄхь. ёс. -- 1954. -- ╥. 34, No 1. -- ╤. 81--88.

\bibitem{16} ╪шЁ°ют └.╚. ╧юфрыухсЁ√ ётюсюфэ√ї ышхт√ї рыухсЁ // ╠рЄхь. ёс. -- 1953. -- ╥. 33, No 2. -- ╤. 441--452.

\bibitem{17} Witt E. Die Unterringe der freien Lieschen Ringe // Math. Z. -- 1956. -- Vol. 64. -- P. 195--216.

\bibitem{18} ╠шїрыхт └.└. ╧юфрыухсЁ√ ётюсюфэ√ї ЎтхЄэ√ї ёєяхЁрыухсЁ ╦ш // ╠рЄ. чрьхЄъш. -- 1985. -- ╥. 37, No 5. -- ╤. 653--661.

\bibitem{19} ╪ЄхЁэ └.╤. ╤тюсюфэ√х ёєяхЁрыухсЁ√ ╦ш // ╤шс. ьрЄ. цєЁэ. -- 1986. -- ╥. 27. -- ╤. 170--174.

\bibitem{Shest1}  Shestakov I.P., Umirbaev U.U. The Nagata automorphism is wild // Proc. Natl. Acad. Sci. USA. -- 2003. -- Vol. 100, No 22. -- P. 12561--12563.

\bibitem{Shest2} Shestakov I.P. and Umirbaev U.U. Tame and wild automorphisms of rings of polynomials in three variables // J. Amer. Math. Soc. -- 2004. -- Vol. 17. -- P. 197--227.

\bibitem{Shest3} ╙ьшЁсрхт ╙.╙., ╪хёЄръют ╚.╧. ╧юфрыухсЁ√ ш ртЄюьюЁЇшчь√ ъюыхЎ ьэюуюўыхэют // ─юъы. ╨└═. -- 2002. -- ╥. 386, No 6. -- ╤. 745--748.

\bibitem{Umir1} ╙ьшЁсрхт ╙.╙. ╬яЁхфхы ■∙хх ёююЄэю°хэш  уЁєяя√ Ёєўэ√ї ртЄюьюЁЇшчьют рыухсЁ√ ьэюуюўыхэют ш фшъшх ртЄюьюЁЇшчь√ ётюсюфэ√ї рёёюЎшрЄштэ√ї рыухсЁ // ─юъы. ╨└═. -- 2006. -- ╥. 407, No 3. -- ╤. 319--324.

\bibitem{Umir2} Umirbaev U.U. The Anick automorphism of free associative algebras // J. Reine Angew. Math. -- 2007. -- Vol. 605. -- P. 165--178.

\bibitem{Kol}  Kolchin E.R. Differential Algebra and Algebraic Groups. Pure and Applied Mathematics, 54. Academic Press,
New York-London, 1973.

\bibitem{KLMP} Kondratieva M.V., Levin A.B., Mikhalev A.V., Pankratiev E.V. Differential and difference dimension polynomials. Mathematics and its Applications, 461. Kluwer Academic Publishers, Dordrecht, 1999.

\bibitem{Ritt}  Ritt J.F. Differential Algebra. Dover Publications, Inc., New York, 1966.

\bibitem{Nag} Nagata M. On Automorphism Group of $k[x,y]$. Lectures in Math., Kyoto Univ., No. 5 (1972) Kinokuniya-Tokyo.

\bibitem{Cohn}  Cohn P.M. Free ideal rings and localization in general rings. New Mathematical Monographs, 3. Cambridge University Press, Cambridge, 2006.

\bibitem{Shir}  ╪шЁ°ют └.╚. ╩юы№Ўр ш рыухсЁ√. ╠.: ═рєър, 1984.

\bibitem{Mag} ╠руэєё ┬., ╩рЁЁрё └., ╤юышЄ¤Ё ─. ╩юьсшэрЄюЁэр  ЄхюЁш  уЁєяя. ╧хЁ. ё рэуы. ╠.: ═рєър, 1974.

\end{thebibliography}
\end{document}